\documentclass[number,dvips]{arxbj}
\usepackage{graphics}

\aid{0}
\volume{16}
\issue{4}
\pubyear{2010}
\firstpage{1177}
\lastpage{1190}
\doi{10.3150/09-BEJ240}

\makeatletter
\newtheorem{Th}{Theorem}[section]
\newtheorem{Co}[Th]{Corollary}
\newtheorem{Lem}[Th]{Lemma}
\newremark{example}{Example}
\makeatother

\begin{document}
\begin{frontmatter}

\title{Asymptotic distributions for a class of generalized $L$-statistics}
\runtitle{Generalized $L$-statistics}

\begin{aug}
\author[1]{\fnms{Yuri V.} \snm{Borovskikh}\thanksref{1}\ead[label=e1]{boryuri@mail.ru}} \and
\author[2]{\fnms{N.C.} \snm{Weber}\corref{}\thanksref{2}\ead[label=e2]{neville.weber@sydney.edu.au}}
\runauthor{Y.V Borovskikh and N.C. Weber}
\address[1]{Department of Applied Mathematics, Transport University,
Moskovsky Avenue 9, 190031, St. Petersburg, Russia. \printead{e1}}
\address[2]{School of Mathematics and Statistics, F07, University of
Sydney, NSW 2006, Australia.\\ \printead{e2}}
\end{aug}

\received{\smonth{9} \syear{2008}}
\revised{\smonth{5} \syear{2009}}

%
\begin{abstract}
We adapt the techniques in Stigler [\textit{Ann. Statist.} \textbf{1}
(1973) 472--477] to obtain a new, general asymptotic
result for trimmed $U$-statistics via the generalized $ L $-statistic
representation introduced by Serfling [\textit{Ann. Statist.} \textbf
{12} (1984) 76--86]. Unlike existing results, we do
not require continuity of an associated distribution at the
truncation points. Our results are quite general
and are expressed in terms of the quantile function associated with the
distribution of the $U$-statistic summands.
This approach leads to improved
conditions for the asymptotic normality of these trimmed $U$-statistics.
\end{abstract}

%

\end{frontmatter}
%

\section{Introduction and statement of results}

Stigler \cite{St} developed an asymptotic result for the trimmed mean without
requiring continuity of the underlying distribution function associated with
the observations. This result was extended to non-degenerate $U$-statistics
based on trimmed samples in Borovskikh and Weber \cite{BW}. An
alternative method
for developing robust versions of $U$-statistics is to consider the statistic
formed by trimming the kernel values, rather than the observations upon
which the
statistic is based. This idea is discussed in, for example, Serfling
\cite{S2},
Choudhury and Serfling \cite{CS} and Gijbels, Janssen and Veraverbeke
\cite{Gi}.
In this paper, we use the generalized $L$-statistic representation developed
in Serfling \cite{S2} to obtain an asymptotic result for trimmed $U$-statistics
under quite general conditions. We will not require continuity of the relevant,
associated distribution at the truncation points.

Let $X,X_1,\ldots, X_n $ be independent identically distributed
random variables, taking values in a measurable space $({\mathsf{X}},  \mathcal{B}(\mathsf{X}))$
and having common distribution $F.$ Let $h$ be a symmetric
function from $ {\mathsf{X}}^m $ to $ R $ and denote by $H_F$ the right-continuous
distribution function of the random variable $h (X_1, \ldots, X_m)$. Set
$N = {n\choose m}$ and let $h_1 , \ldots,h_N $
be an enumeration of the values of $h(X_{i_{1}}, \ldots, X_{i_{m}})$ taken
over the $ N $ $m$-tuples in $\sigma_{nm} = \{ (i_1 , \ldots, i_m)\dvtx 1 \leq i_1 < \cdots< i_m \leq n \}.$ Note that these random
variables $ h_i$
are, in general, dependent. Let $h_{n1} \leq\cdots\leq h_{nN}$ denote the
ordered values of $ h_1 , \ldots, h_N $.

The original $U$-statistic is defined as an average taken over the $N$
possible outcomes $h (X_{i_{1}}, \ldots, X_{i_{m}}),  1 \leq i_1
<\cdots< i_m \leq n$,
that is,
%
\begin{equation}\label{for1}
U =  \pmatrix{ n\cr m}^{-1} \sum_{\sigma_{nm}} h (X_{i_{1}}, \ldots,
X_{i_{m}}) =  N^{-1} \sum_{i=1}^N h_{ni}
= \int_{R} x\, \mathrm{d}H_n(x) ,
\end{equation}
where the empirical distribution function $ H_n(x)$ of $U$-statistical
structure is defined by
%
\begin{equation}\label{for2}
H_n (y) = \pmatrix{ n\cr m}^{-1} \sum_{\sigma_{nm}}
I \{ h ( X_{i_1}, \ldots, X_{i_m}) \leq y\},\qquad    y \in R,
\end{equation}
and $I\{ A \}$ denotes the indicator of the set $A$.
For any $0<\gamma<1$, let $ N_\gamma= [\gamma N]$,
where $[  a  ]$ denotes the largest integer less than or equal to $ a.$
If $0 < \alpha< \beta< 1$,
then put $N_{\alpha\beta} = N_\beta- N_\alpha .$ The trimmed
versions of
$U$ are based on trimming the second sum in (1),
%
\begin{equation}\label{for3}
{ U_{\alpha\beta} } = N_{\alpha\beta}^{-1} \sum_{i = N_{\alpha} +
1}^{N_{\beta}} h_{ni} ,
\end{equation}
or on trimming of the range of integration in (\ref{for1}),
%
\begin{equation}\label{for4}
L_{\alpha\beta} = \int_{[h_\alpha, h_\beta)} x \, \mathrm{d}H_n(x) ,
\end{equation}
with $h_\alpha= h_{n\bar N_\alpha}$ and $ h_\beta= h_{n\bar N_\beta},$
where $\bar N_\gamma= -[-\gamma N] ,\gamma= \alpha,\beta.$ For the results
that follow, it is important to note that the lower
bound for the integral in (\ref{for4}) is included and the upper bound excluded.
This is
critical since $H_n$ is a step function. With this constraint, we are
able to obtain
the asymptotic distribution of ${L_{\alpha\beta}}$ without imposing
any conditions on the nature
of $H_F.$ In Lemma \ref{lem2.3}, we show that
$ L_{\alpha\beta} = N^{-1}\sum_{i = \bar{N}_{\alpha}}^{\bar
{N}_{\beta} -1}  h_{ni}.$ Thus, $U_{\alpha\beta}$ and $L_{\alpha\beta}$ differ in
terms of their
divisors, and there are possible subtle differences in the number of summands.


A class of generalized $L$-statistics, which includes (\ref{for3}) and (\ref{for4}), was
introduced by
Serfling~\cite{S2}. The trimmed $U$-statistics (\ref{for3}) and (\ref{for4}) are
directly connected
with generalized Lorenz curves, which are important in financial
mathematics (see, for example, Goldie \cite{Go}, Helmers and Zitikis~\cite{HZ}).

Clearly, $H_n (y)$ is an unbiased estimator of $H_F (y).$ In the case
$m=1$ and $h (x) = x,$ $H_n$ reduces to the usual empirical distribution
function.
Define the left-continuous quantile function $H^{-1}_F (t) =
\inf\{y \in R \dvtx  H_F (y) \geq t \},  0 <t \leq1,   H^{-1}_F(0)=H^{-1}_F(0+),$
for any distribution function~$H_F.$
The empirical quantile function $ H^{-1}_n(t) $ has the form
\[
H^{-1}_n(t) = \sum_{i=1}^N h_{ni}I\biggl\{\frac{i-1}{N}< t \leq\frac
{i}{N}\biggr\}, \qquad   H^{-1}_n(0) = h_{n1}.
\]
A large number of authors have studied the weak convergence of such
$L$-statistics in the case $ m=1, h(x)=x.$
A partial list consists of Chernoff \textit{et al.}~\cite{CG}, Bickel
\cite{Bi},
Shorack \cite{Sh,Sh1}, Stigler \cite{St,St1}, Cs\"{o}rgo \textit{et al.}~\cite{Cs},
Griffin and Pruitt \cite{Gr}, Cheng \cite{Ch}, Mason and Shorack
\cite{Ma}.
For $ m \geq2 $, under various sets of regularity conditions, asymptotic
normality of various types of generalized $L$-statistics has been investigated
by Silverman \cite{Si}, Serfling \cite{S2}, Akritas \cite{Ak},
Janssen \textit{et
al.}~\cite{JS}, Helmers
and Ruymgaart \cite{HR}, Gijbels \textit{et al.}~\cite{Gi} and H\"
{o}ssjer \cite{Ho}.

In the aforementioned papers, for $ m\geq2 $, the results always
assumed that
$ H_F $ is continuous or smooth. However, in modern statistical robust
procedures
and for bootstrap procedures, results allowing for the discontinuity of the
underlying distribution function $ H_F $ are needed. We study the asymptotic
behavior of $ U_{\alpha\beta} $ and $ L_{\alpha\beta}$ for any $
H_F $ without imposing the requirement of
continuity.

The conditions of our theorem and the limit random variable are defined via
the values of quantile function $ H^{-1}_F $ at the points $ \alpha$ and
$ \beta.$ Existing results handle the cases where
$ H^{-1}_F(\gamma+) = H^{-1}_F(\gamma) , \gamma= \alpha,\beta.$ Our
main result is derived without this assumption of continuity. We represent
the trimmed $U$-statistic as a sum of classical $U$-statistics with bounded,
non-degenerate kernels plus some smaller terms and then we apply standard
results to such statistics.

For convenience, in what follows, for the distribution function $ H_F,
$ we denote
the smallest quantile $ H^{-1}_F(\gamma) $ and the largest quantile
$ H^{-1}_F(\gamma+)$ as, respectively,
\begin{eqnarray*}
\xi^-_\gamma&:=& \inf\{x\in R \dvtx  H_F(x)\geq\gamma\},\\
  \xi
^+_\gamma&:=& \sup\{x\in R \dvtx  H_F(x)\leq\gamma\}
\end{eqnarray*}
and $ \Delta\xi_\gamma= \xi^+_\gamma- \xi^-_\gamma$ with
$ \gamma = \alpha, \beta$. Let
\[
\dot N^{\pm}_\gamma = \sum_{i=1}^N I\{h_i < \xi^{\pm}_\gamma\},
\qquad
N^{\pm}_\gamma = \sum_{i=1}^N I\{h_i \leq\xi^{\pm}_\gamma\}.
\]
Note that
%
\begin{equation}\label{for5}
H_n(\xi^\pm_\gamma-) = N^{-1}\dot N^\pm_\gamma,\qquad
H_n(\xi^\pm_\gamma) = N^{-1}N^\pm_\gamma
\end{equation}
and $ H^{-1}_n(\gamma)=h_{n\bar N_\gamma} $ are valid for all $
0<\gamma<1 $
and the following events coincide:
\begin{equation}\label{for6}
\{H^{-1}_n(\gamma) > x\} = \{\gamma> H_n(x)\}, \qquad   \{H^{-1}_n(\gamma
) \leq x\} =
\{\gamma\leq H_n(x)\},\qquad   x \in R.\nonumber
\end{equation}
Introduce the functional $\theta= \theta(H_F)$, where
\[
\theta= \int_{R}\bigl[\bigl((x-\xi^-_\beta)I\{x\leq\xi^-_\beta\} +\beta
\xi^-_\beta\bigr) - \bigl((x-\xi^+_\alpha)I\{x < \xi^+_\alpha\} +\alpha\xi
^+_\alpha\bigr) \bigr]\,\mathrm{d}H_F(x)
\]
and the following functions with $ x \in\mathsf{X} $:
%
\begin{eqnarray}\label{for7}
g(x) & = &\bigl[ EI\{h(x,X_2,\ldots,X_m) \leq\xi^-_\beta\}
\bigl(h(x,X_2,\ldots,X_m) - \xi^-_\beta\bigr)+\beta\xi^-_\beta\bigr] \nonumber\\
& &{} - \bigl[ EI\{ h(x,X_2,\ldots,X_m) < \xi^+_\alpha\}\bigl(h(x,X_2,\ldots
,X_m) - \xi^+_\alpha\bigr)+\alpha\xi^+_\alpha\bigr] -\theta,\nonumber\\
g_\alpha(x) & = & EI\{h(x,X_2,\ldots,X_m) < \xi^+_\alpha\}-\theta
_\alpha, \qquad  \theta_\alpha=H_F(\xi^+_\alpha-),\\
g_\beta(x) & = & EI\{h(x,X_2,\ldots,X_m) \leq\xi^-_\beta\}-\theta
_\beta \nonumber\\
& = & 1-\theta_\beta- EI\{h(x,X_2,\ldots,X_m) > \xi^-_\beta\},\qquad
 \theta_\beta=H_F(\xi^-_\beta).\nonumber
\end{eqnarray}
Note that for all $ 0<\alpha<\beta<1 $ and $ x \in\mathsf{X} $, we have
$ |g(x)| \leq4(|\xi^+_\alpha|+|\xi^-_\beta|).$

Let $   \sigma_g^2 = Eg^2(X),   \sigma_{g_\alpha}^2 = Eg_\alpha^2(X),
  \sigma_{g_\beta}^2 = Eg_\beta^2(X),
c_{gg_{\alpha}} = Eg(X)g_{\alpha}(X),
c_{gg_{\beta}} = Eg(X)g_{\beta}(X)$   and
$c_{g_{\alpha}g_{\beta}} = Eg_{\alpha}(X)g_{\beta}(X).
$

\begin{Th}\label{teo1.1}
If $ \sigma_{g}^2 > 0,$ then for any underlying distribution function
$ H_F$, we
have
\[
\frac{\beta-\alpha}{m}\sqrt{n}( U_{\alpha\beta} - \theta)
\stackrel{d}{\longrightarrow} \tau_{g}
- \Delta\xi_\alpha  I(\tau_\alpha> 0) \tau_\alpha
- \Delta\xi_\beta  I(\tau_\beta< 0) \tau_\beta,
\]
where $ (\tau_\alpha, \tau_g ,\tau_\beta) $ is a trivariate
Gaussian random vector
with mean vector zero and covariance matrix
\[
\pmatrix{
\sigma_{g_\alpha}^2 & c_{gg_{\alpha}} & c_{g_{\alpha}g_{\beta}} \vspace*{2pt}\cr
c_{gg_{\alpha}} & \sigma_g^2 & c_{gg_{\beta}} \vspace*{2pt}\cr
c_{g_{\alpha}g_{\beta}} & c_{gg_{\beta}} & \sigma_{g_\beta}^2}.
\]
\end{Th}

\begin{Co}\label{cor1.2}
For any underlying distribution function $ H_F $, we have, when $
\sigma_{g}^2 > 0 $,
\[
\frac{\sqrt{n}}{m}( L_{\alpha\beta} -\theta) \stackrel
{d}{\longrightarrow} \tau_{g}
- \Delta\xi_\alpha  I(\tau_\alpha> 0) \tau_\alpha
- \Delta\xi_\beta  I(\tau_\beta< 0) \tau_\beta,
\]
where $ (\tau_\alpha, \tau_g ,\tau_\beta) $ is a trivariate
Gaussian random vector defined as in Theorem \ref{teo1.1}.
\end{Co}

\begin{Co}\label{cor1.3}
Suppose that the quantile function $ H^{-1}_{F}(x) $
is continuous at the points $ \alpha$ and $ \beta.$
If $ \sigma^2_g >0,$ then
\[
\frac{\beta-\alpha}{m}\sqrt{n}( U_{\alpha\beta} - \theta)
\stackrel{d}{\longrightarrow}
\tau_g .
\]
\end{Co}

For the simple case $ m = 1 $, the functions in (\ref{for7}) reduce to
\begin{eqnarray*}
g(x) & = & I\{\xi^+_\alpha\leq h(x) \leq\xi^-_\beta\}h(x) - EI\{
\xi^+_\alpha\leq h(X) \leq\xi^-_\beta\}h(X) \nonumber\\
& & {}+ \xi^+_\alpha g_\alpha(x) - \xi^-_\beta g_\beta(x) ,\nonumber
\\
g_\alpha(x) & = & I\{h(x) < \xi^+_\alpha\}-\theta_\alpha,
\nonumber\\
g_\beta(x) & = & I\{h(x) \leq\xi^-_\beta\}-\theta_\beta
= 1-\theta_\beta- I\{h(x) > \xi^-_\beta\}.
\end{eqnarray*}
A useful application of the theorem for the $m=2$ case is for the
kernel $ h(x,y) $ $ = \frac{1}{2}(x-y)^2 . $
This provides the asymptotic behavior of a natural, alternative robust version
of the sample variance. We will now develop explicit expressions for
the terms in a more
interesting example.

\begin{example*}
 Let $ h(x_1,\ldots,x_m)=\max\{x_1,\ldots
,x_m\} $ with $ m \geq2 $ . Let $ F(t) $ be the distribution function
of $ X_i $ and
let $ Y=\max\{X_2,\ldots,X_m \}.$
Then $ H_F(t) = (F(t))^m $ and
\begin{eqnarray*}
g(x) & = & g_{\alpha\beta}(x) + \xi^+_\alpha g_\alpha(x) - \xi
^-_\beta g_\beta(x)  , \nonumber\\
g_{\alpha\beta}(x) & = & EI\bigl\{\xi^+_\alpha\leq\max\{x,Y\} \leq\xi
^-_\beta\bigr\}\max\{x,Y\} -
\int_{[\xi^+_\alpha,\xi^-_\beta]}y\,\mathrm{d}H_F(y) \nonumber\\
& = & I\{\xi^+_\alpha\leq x \leq\xi^-_\beta\} x (F(x))^{m-1} -
\int_{[\xi^+_\alpha,\xi^-_\beta]} y (F(y))^{m-1} \,\mathrm{d} F(y) \nonumber
\\
& & {}+ \int_{[\xi^+_\alpha,\xi^-_\beta]} \bigl(I\{x<y\} - F(y-)\bigr)y
\,\mathrm{d}(F(y))^{m-1} , \nonumber\\
g_\alpha(x) & = & I\{x < \xi^+_\alpha\}(F(\xi^+_\alpha-))^{m-1} -
(F(\xi^+_\alpha-))^{m},  \nonumber\\
g_\beta(x) & = & I\{x \leq
\xi^-_\beta\}(F(\xi^-_\beta))^{m-1}-(F(\xi^-_\beta))^{m}.
\end{eqnarray*}
In addition,
\begin{eqnarray*}
\sigma^2_g &=& Eg_{\alpha\beta}^2(X) + E[\xi^+_\alpha g_\alpha(X) -
\xi^-_\beta g_\beta(X)]^2\\
&&{} + 2Eg_{\alpha\beta}(X) [\xi^+_\alpha
g_\alpha(X) -
\xi^-_\beta g_\beta(X)],\\
\sigma_{g_\alpha}^2 &=& (F(\xi^+_\alpha-))^{2m-1}\bigl(1 - F(\xi^+_\alpha-)\bigr),\qquad
 \sigma_{g_\beta}^2 = (F(\xi^-_\beta))^{2m-1}\bigl(1 - F(\xi
^-_\beta)\bigr) ,\\
Eg_{\alpha}(X)g_{\beta}(X) &=& (F(\xi^+_\alpha-)F(\xi^-_\beta
))^{m-1} F(\xi^+_\alpha-) \bigl(1 - F(\xi^-_\beta)\bigr)  ,\\
Eg_{\alpha\beta}(X)g_{\alpha}(X) & = & (F(\xi^+_\alpha-))^{m} \int
_{[\xi^+_\alpha,\xi^-_\beta]} \bigl(1 - F(y-)\bigr)y\, \mathrm{d}(F(y))^{m-1} \nonumber
\\
& & {}- (F(\xi^+_\alpha-))^{m} \int_{[\xi^+_\alpha,\xi^-_\beta]} y
(F(y))^{m-1}\, \mathrm{d} F(y),\\
Eg_{\alpha\beta}(X)g_{\beta}(X) & = & (F(\xi^-_\beta
))^{m-1}\bigl(1-F(\xi^-_\beta)\bigr)
\int_{[\xi^+_\alpha,\xi^-_\beta]} y (F(y))^{m-1}\, \mathrm{d} F(y) \nonumber
\\
& &{} + (F(\xi^-_\beta))^{m-1}\bigl(1-F(\xi^-_\beta)\bigr)
\int_{[\xi^+_\alpha,\xi^-_\beta]} F(y-) y\, \mathrm{d}(F(y))^{m-1}.
\end{eqnarray*}

Consider the distribution function
\begin{eqnarray*}
F(t) & = & 2tI\bigl\{0\leq t < \tfrac{1}{2}\alpha^{1/m}\bigr\}+ \alpha^{1/m}
I\bigl\{ \tfrac{1}{2} \alpha^{1/m} \leq t < \alpha^{1/m}\bigr\} \nonumber\\
& &{} + tI\{\alpha^{1/m} \leq t < \beta^{1/m} \} + \beta^{1/m}I\{
\beta^{1/m} \leq t < 2 \beta^{1/m}\} \nonumber\\
& &{} + \tfrac{1}{2}tI\{2 \beta^{1/m} \leq
t < 2 \} + I\{ t \geq2 \} , \qquad   t \in R  .
\end{eqnarray*}
Then
\begin{eqnarray*}
\xi^-_\alpha&=& \tfrac{1}{2} \alpha^{1/m},  \qquad  \xi^+_\alpha=
\alpha^{1/m},\qquad
 \xi^-_\beta= \beta^{1/m},\qquad   \xi^+_\beta= 2\beta^{1/m},\\
F(\xi^+_\alpha-)&=& \alpha^{1/m} ,\qquad   F(\xi^-_\beta) = \beta
^{1/m} ,\qquad    F(t)=t ,\qquad   t \in[\alpha^{1/m},\beta^{1/m}] ,\qquad
\sigma^2_g > 0
\end{eqnarray*}
and the limiting behavior is given by
\[
\frac{\beta-\alpha}{m}\sqrt{n}( U_{\alpha\beta} - \theta)
\stackrel{d}{\longrightarrow} \tau_{g}
- \frac{1}{2} \alpha^{1/m}   I(\tau_\alpha> 0) \tau_\alpha
- \beta^{1/m}   I(\tau_\beta< 0) \tau_\beta .
\]

However, for the simpler distribution function
\[
F(t) = tI\{0 \leq t < 1 \} + I\{t\geq1\} ,\qquad   t \in R  ,
\]
we have
\begin{eqnarray*}
\xi^-_\alpha&=& \xi^+_\alpha= \alpha^{1/m},\qquad
 \xi^-_\beta= \xi^+_\beta= \beta^{1/m},\\
F(\xi^+_\alpha-) &=& \alpha^{1/m},\qquad   F(\xi^-_\beta) = \beta
^{1/m},\qquad    F(t)=t, \qquad  t \in[\alpha^{1/m},\beta^{1/m}] ,\qquad
\sigma^2_g > 0
\end{eqnarray*}
and we get the asymptotic behavior covered by Janssen \textit{et
al.}~\cite{JS},
\[
\frac{\beta-\alpha}{m}\sqrt{n}( U_{\alpha\beta} - \theta)
\stackrel{d}{\longrightarrow} \tau_{g} .
\]
\end{example*}

\section{Proofs}
The following two lemmas are key results for the proof.

\begin{Lem}\label{lem2.1}
The following representation holds:
%
\begin{eqnarray}\label{for8}
\sum_{i = N_{\alpha} + 1}^{N_{\beta}} h_{ni} & = & \sum_{i = 1}^{N}
I\{\xi^+_\alpha\leq h_i \leq\xi^-_\beta\}h_i +
\xi^+_\alpha(\dot N^+_\alpha-N_\alpha)- \xi^-_\beta(N^-_\beta
-N_\beta) \nonumber\\
& &{} -\Delta\xi_\alpha I\{N_\alpha< \dot N^+_\alpha\}(\dot
N^+_\alpha-N_\alpha)-
\Delta\xi_\beta I\{N^-_\beta< N_\beta\}(N^-_\beta-N_\beta)\nonumber
\\
& &{} +\mathbb{L}_\alpha+ \mathbb{L}_\beta ,
\end{eqnarray}
where $ \mathbb{L}_\alpha= J_\alpha- \bar J_\alpha$ with
\[
J_\alpha= I\{N_\alpha< \dot N^+_\alpha\}\sum_{i=N_\alpha+1}^{\dot
N^+_\alpha}(h_{ni}-\xi^-_\alpha) ,\qquad
\bar J_\alpha= I\{\dot N^+_\alpha\leq N_\alpha\}\sum_{i=\dot
N^+_\alpha+1}^{N_\alpha}(h_{ni}-\xi^+_\alpha)
\]
and $ \mathbb{L}_\beta = \bar J_\beta- J_\beta$ with
\begin{eqnarray*}
J_\beta&=& I\{N_\beta< N^-_\beta\}\sum_{i=N_\beta+1}^{N^-_\beta
}(h_{ni}-\xi^-_\beta) ,\\
\bar J_\beta&=& I\{N^-_\beta\leq N_\beta\}\sum_{i=N^-_\beta
+1}^{N_\beta}(h_{ni}-\xi^+_\beta) .
\end{eqnarray*}
\end{Lem}

\begin{pf}
 For $ i = 1,\ldots,N $, write
\[
\dot h_{ni} = (h_{ni}+\Delta\xi_\alpha)I\{h_{ni} < \xi^+_\alpha\}+
h_{ni}I\{\xi^+_\alpha\leq h_{ni} \leq\xi^-_\beta\}+
(h_{ni}-\Delta\xi_\beta)I\{h_{ni} > \xi^-_\beta\} .
\]
Since $ \dot h_{ni} = h_{ni}+\Delta\xi_\alpha I\{h_{ni} < \xi
^+_\alpha\}-\Delta\xi_\beta I\{h_{ni} > \xi^-_\beta\} $,
$ I\{h_{ni} < \xi^+_\alpha\} $ = $ I\{i \leq\dot N^+_\alpha\} $
and, by (\ref{for6}),
$ I\{h_{ni} > \xi^-_\beta\} = I\{i > N^-_\beta\},$ we can write
\begin{eqnarray}\label{for9}
\sum_{i = N_{\alpha} + 1}^{N_{\beta}} h_{ni} & = & \sum_{i =
N_{\alpha} + 1}^{N_{\beta}}\dot h_{ni}
- \Delta\xi_\alpha I\{N_\alpha< \dot N^+_\alpha\}(\dot N^+_\alpha
-N_\alpha)\nonumber \\ 
& & {}- \Delta\xi_\beta I\{N^-_\beta< N_\beta\}(N^-_\beta-N_\beta) .
\end{eqnarray}
Note that $ h_{n\dot N^+_\alpha} < \xi^+_\alpha\leq h_{n,\dot
N^+_\alpha+1} $
and $ h_{nN^-_\beta} \leq\xi^-_\beta< h_{n,N^-_\beta+1}.$ From
(\ref{for6}), we have
$I\{\xi^+_\alpha\leq h_{ni} \leq\xi^-_\beta\} = $ $ I\{\dot
N^+_\alpha < i \leq
N^-_\beta\}.$ Hence, in (\ref{for9}),
%
\begin{eqnarray}\label{for10}
\sum_{i = N_{\alpha} + 1}^{N_{\beta}} \dot h_{ni} & = & \sum_{i =
\dot N^+_\alpha+ 1}^{N^-_{\beta}} h_{ni}
- I\{\dot N^+_\alpha\leq N_\alpha\}\sum_{i=\dot N^+_\alpha
+1}^{N_\alpha}h_{ni} \nonumber\\
& &{} +I\{N_\alpha< \dot N^+_\alpha\}\sum_{i=N_\alpha+1}^{\dot
N^+_\alpha}(h_{ni}+\Delta\xi_\alpha)
- I\{N_\beta< N^-_\beta\}\sum_{i=N_\beta+1}^{N^-_\beta}h_{ni}
\nonumber\\
& &{} + I\{N^-_\beta\leq N_\beta\}\sum_{i=N^-_\beta+1}^{N_\beta
}(h_{ni}-\Delta\xi_\beta)\nonumber \\
& = & \sum_{i = 1}^{N} I\{\xi^+_\alpha\leq h_i \leq\xi^-_\beta\}h_i
+ \xi^+_\alpha(\dot N^+_\alpha-N_\alpha)- \xi^-_\beta(N^-_\beta
-N_\beta) \nonumber\\
& &{} + \mathbb{L}_\alpha+ \mathbb{L}_\beta .
\end{eqnarray}
Equation (\ref{for8}) follows from (\ref{for9}) and (\ref{for10}). This proves Lemma \ref{lem2.1}.
\end{pf}

\begin{Lem}\label{lem2.2}
Note that
\begin{eqnarray*}
N^{-1}\sum_{i = N_{\alpha} + 1}^{N_{\beta}} h_{ni} & = &
N^{-1} \sum_{i = 1}^{N} I\{\xi^+_\alpha\leq h_i \leq\xi^-_\beta\}
h_i +
\xi^+_\alpha\bigl(H_n(\xi^+_\alpha-) - \alpha\bigr) \nonumber\\
& &{} - \xi^-_\beta\bigl(H_n(\xi^-_\beta) - \beta\bigr)
- \Delta\xi_\alpha I\{N_\alpha< \dot N^+_\alpha\}
\bigl(H_n(\xi^+_\alpha-) -
\alpha\bigr) \nonumber\\
& &{} - \Delta\xi_\beta I\{N^-_\beta< N_\beta\}  \bigl(H_n(\xi
^-_\beta) - \beta
\bigr)
+ n^{-1/2} \varrho_n ,
\end{eqnarray*}
where $ \varrho_n \rightarrow0 $
in probability as $ n \rightarrow\infty .$
\end{Lem}

\begin{pf}
 We shall estimate $ \mathbb
{L}_\alpha$ and
$ \mathbb{L}_\beta,$ taking into account the values of the distribution
function~$ H_F(x) $ at $ x=\xi^{\pm}_\gamma$ with $ \gamma=\alpha
,\beta.$
Figures \ref{fig1} and \ref{fig2} illustrate the different situations that need to be
considered.

\begin{figure}

\includegraphics{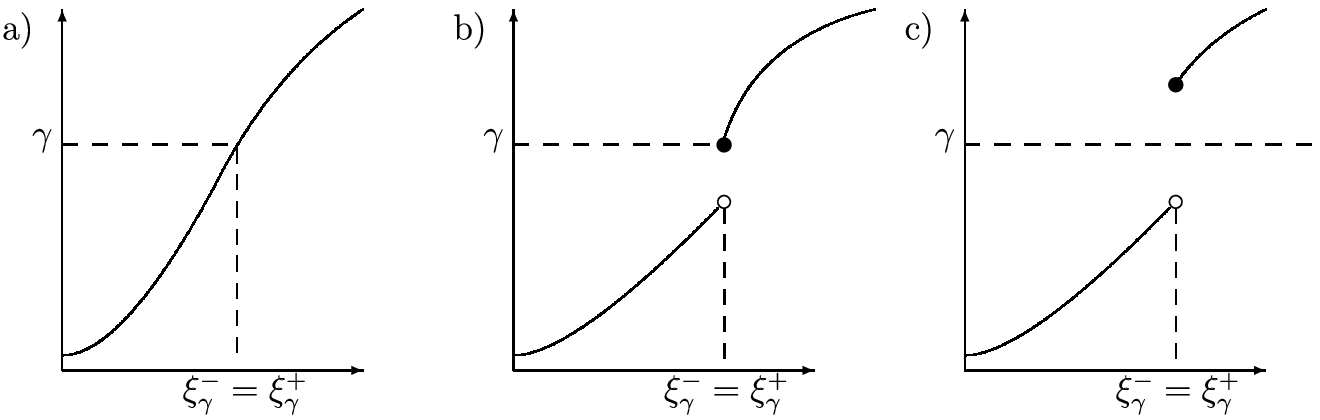}

  \caption{Plots of $H(\cdot)$ with $\xi_{\gamma}^{-} = \xi_{\gamma}^{+}$:
(a) $H(\xi_{\gamma}^{\pm}-) = \gamma=H(\xi_{\gamma}^{\pm})$;
(b) $H(\xi_{\gamma}^{\pm}-) < \gamma=H(\xi_{\gamma}^{\pm})$;
 (c)~$H(\xi_{\gamma}^{\pm}-) < \gamma<H(\xi_{\gamma}^{\pm}).$}\label{fig1}
\end{figure}

\begin{figure}[b]

\includegraphics{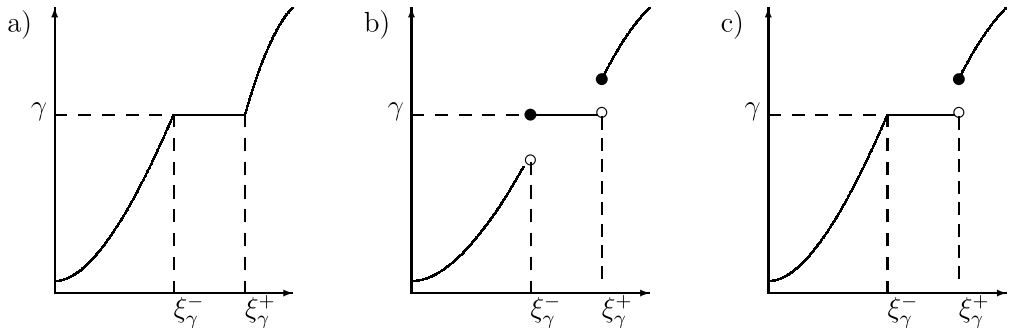}

  \caption{Plots of $H(\cdot)$ with $\xi_{\gamma}^{-} < \xi_{\gamma}^{+}$:
(a) $\gamma=H(\xi_{\gamma}^{-}) = H(\xi_{\gamma}^{+})$;  (b)
$H(\xi_{\gamma}^{-} -) < \gamma=  H(\xi_{\gamma}^{-}) =\break H(\xi
_{\gamma}^{+}-) < H(\xi_{\gamma}^{+})$; (c)
$H(\xi_{\gamma}^{-}) = \gamma= H(\xi_{\gamma}^{+}-) < H(\xi
_{\gamma}^{+}).$}\label{fig2}
\end{figure}

%
%
%
%
%
%
%
%
%

\textbf{Estimating $ \mathbb{L}_\alpha$}. Noting that $ I\{\xi^+_\alpha
> h_{ni} \} = I\{i\leq\dot N^+_\alpha\} $,
we write
\begin{eqnarray*}
J_\alpha& = & I\{N_\alpha< \dot N^+_\alpha\}\sum_{i=N_\alpha
+1}^{\dot N^+_\alpha}(h_{ni}-\xi^-_\alpha)I\{\xi^+_\alpha> h_{ni}
\} \nonumber\\
& = & I\{N_\alpha< \dot N^+_\alpha\}\sum_{i=N_\alpha+1}^{\dot
N^+_\alpha}(h_{ni}-
\xi^-_\alpha)I\{\xi^-_\alpha< h_{ni} < \xi^+_\alpha\} I\{
N^-_\alpha< i \leq\dot N^+_\alpha\}\nonumber\\
& & {}- I\{N_\alpha< \dot N^+_\alpha\}\sum_{i=N_\alpha+1}^{\dot
N^+_\alpha}(
\xi^-_\alpha- h_{ni})I\{\xi^-_\alpha\geq h_{ni} \} I\{ i \leq
N^-_\alpha\}\\
&=& J^+_\alpha- J^-_\alpha.
\end{eqnarray*}
It is clear that if $ \xi^-_\alpha= \xi^+_\alpha$, then $
J^+_\alpha= 0 $ a.s.
Let $ \xi^-_\alpha\neq\xi^+_\alpha$, as is the case in Fig. \ref{fig2}.
In this case,
$ H_F(\xi^-_\alpha) =  \alpha =  H_F(\xi^+_\alpha-) $ and we
can write
%
\begin{eqnarray*}
0 &\leq& J^+_\alpha \leq I\{N^-_\alpha< \dot N^+_\alpha\}(\xi
^+_\alpha- \xi^-_\alpha)\sum_{i=N^-_\alpha+1}^{\dot N^+_\alpha}
I\{\xi^-_\alpha< h_{ni} < \xi^+_\alpha\} \\
& \leq& \Delta\xi_\alpha \sum_{i=1}^{N}
I\{\xi^-_\alpha< h_i < \xi^+_\alpha\} = 0\qquad     \mbox{a.s.}
\end{eqnarray*}
since $ EI\{\xi^-_\alpha< h_i < \xi^+_\alpha\} = H_F(\xi^+_\alpha
-) -
H_F(\xi^-_\alpha) = 0.$ Hence, we always have the
relation\vspace*{1pt} \mbox{$ J^+_\alpha= 0 $} a.s. To estimate $ J^-_\alpha$, we note
first that if $ \dot N^+_\alpha> N^-_\alpha$, then
the indicator $ I\{ i \leq N^-_\alpha\} = 0 $ for all $ i = N^-_\alpha
+1, \ldots,
\dot N^+_\alpha.$ Therefore, we have the inequalities
\begin{eqnarray}\label{for11}
0 &\leq& J^-_\alpha \leq  I\{N_\alpha< N^-_\alpha\}\sum
_{i=N_\alpha+1}^{N^-_\alpha}(
\xi^-_\alpha- h_{ni})I\{\xi^-_\alpha\geq h_{ni} \} \nonumber\\ 
& \leq& I\{N_\alpha< N^-_\alpha\}(N^-_\alpha - N_\alpha)
(\xi^-_\alpha- h_{nN_\alpha})I\{\xi^-_\alpha\geq h_{nN_\alpha} \}.
\end{eqnarray}
Further, we shall apply the technique used in Smirnov \cite{Sm}
with a probability inequality from Hoeffding \cite{Ho} (or see, for example,
Serfling \cite{S1}, pages 75 and 201). Thus,
%
\begin{eqnarray}\label{for12}
&&P\{ (\xi^-_\alpha- h_{nN_\alpha} > \varepsilon) \cap(
\xi^-_\alpha\geq
h_{nN_\alpha} )\} \nonumber\\
&&\quad  \leq  P\{ \xi^-_\alpha - h_{nN_\alpha} \geq\varepsilon\}
\nonumber\\
&&\quad  =  P\{ H_n(\xi^-_\alpha-\varepsilon) \geq N^{-1}N_\alpha\}
\nonumber\\
&&\quad  =  P\{ H_n(\xi^-_\alpha-\varepsilon) - H(\xi^-_\alpha
-\varepsilon) \geq N^{-1}N_\alpha- H(\xi^-_\alpha-\varepsilon)\}
\nonumber\\
&&\quad  \leq c_1 \exp\{ -c_2 n \theta^2_\alpha(\xi^-_\alpha
,\varepsilon)\}
\end{eqnarray}
with some positive constants $ c_1 $ and $ c_2,$ depending only on $ m
$ and
$ \theta_\alpha(\xi^-_\alpha,\varepsilon) = \alpha-\break  H_F(\xi
^-_\alpha-
\varepsilon).$ Further, $ \theta_\alpha(\xi^-_\alpha, $ $
\varepsilon) > 0 $
for any small values of $ \varepsilon> 0 $, by the definition of the
smallest $\alpha$-quantile
$ \xi^-_\alpha.$
Under the conditions of the lemma,
$ \sqrt{n}N^{-1}(N^-_\alpha- N_\alpha) \stackrel{d}{\longrightarrow}
\tau^-_\alpha$ as $ n \rightarrow\infty.$ Hence,
$ \sqrt{n}N^{-1} J_\alpha\rightarrow0 $ in probability as $ n
\rightarrow
\infty.$

Next, we consider $ \bar J_\alpha. $
By definition, $ \dot N^+_\alpha\leq N^+_\alpha$ and
since $ I\{h_{ni}< \xi^+_\alpha\} = I\{i\leq\dot N^+_\alpha\} $ and
$ h_{nN^+_\alpha} \leq\xi^+_\alpha< h_{n,N^+_\alpha+1} $, it
follows that
the indicator $ I\{h_{ni} = \xi^+_\alpha \} = 1 $ for
$ i = \dot N^+_\alpha+1,\ldots, N^+_\alpha$ and we can write
%
\begin{eqnarray}\label{for13}
0 &\leq&\bar J_\alpha =  I\{\dot N^+_\alpha\leq N_\alpha\}\sum
_{i=\dot N^+_\alpha+1}^{ N_\alpha}(h_{ni}-\xi^+_\alpha)
I\{h_{ni} \geq\xi^+_\alpha \} \nonumber\\
& = & I\{ N^+_\alpha\leq N_\alpha\}\sum_{i= N^+_\alpha+1}^{
N_\alpha}(h_{ni}-\xi^+_\alpha)
I\{h_{ni} > \xi^+_\alpha \} \nonumber\\
& \leq & I\{ N^+_\alpha\leq N_\alpha\}(N_\alpha- N^+_\alpha
)(h_{nN_\alpha}-\xi^+_\alpha)
I\{h_{nN_\alpha} > \xi^+_\alpha \} .
\end{eqnarray}
In (\ref{for13}), we need to consider two cases: $ H_F(\xi^+_\alpha) = \alpha
$ and $
\alpha< H_F(\xi^+_\alpha).$ In the first case,
$ H_F(\xi^+_\alpha) = \alpha$ and we have the weak convergence
$ \sqrt{n}N^{-1}(N^+_\alpha- N_\alpha) \stackrel{d}{\longrightarrow
} \tau^+_\alpha$ as $ n \rightarrow\infty $
and the following estimates which are similar to (\ref{for12}):
%
\begin{eqnarray}\label{for14}
&&P\{ (h_{nN_\alpha}-\xi^+_\alpha> \varepsilon) \cap
(h_{nN_\alpha} >\xi^-_\alpha)\}  \nonumber\\
&&\quad  \leq  P\{ h_{nN_\alpha} - \xi^+_\alpha> \varepsilon\} \nonumber
\\
&&\quad =  P\{ N^{-1}N_\alpha> H_n(\xi^+_\alpha+\varepsilon) \}\nonumber
\\
&&\quad  =  P\{ H(\xi^+_\alpha+\varepsilon) - H_n(\xi^+_\alpha
+\varepsilon) > H(\xi^+_\alpha+\varepsilon)-N^{-1}N_\alpha\}
\nonumber\\
&&\quad  \leq c_1 \exp\{ -c_2 n \delta^2_\alpha(\xi^+_\alpha
,\varepsilon)\},
\end{eqnarray}
where $ \delta_\alpha(\xi^+_\alpha, \varepsilon) = H_F(\xi
^+_\alpha+
\varepsilon) - \alpha.$ In addition,
$ \delta_\alpha(\xi^+_\alpha, \varepsilon) > 0 $
for any small values of $ \varepsilon> 0 $ because of the definition
of the
largest $ \alpha$-quantile $ \xi^+_\alpha.$
Hence, in the case $ H_F(\xi^+_\alpha) = \alpha$, we have
$ \sqrt{n}N^{-1} \bar J_\alpha\rightarrow0$ in probability as $ n
\rightarrow
\infty.$
In the second case in (\ref{for13}), $ \delta_\alpha(\xi^+_\alpha,0) =
H_F(\xi^+_\alpha)
- \alpha > 0 $ and we have the representation
%
\begin{equation}\label{for15}
\sqrt{n}N^{-1}(N^+_\alpha- N_\alpha) = \sqrt{n}\bigl(H_n(\xi^+_\alpha
)-H_F(\xi^+_\alpha)\bigr)+\sqrt{n}\delta_\alpha(\xi^+_\alpha,0)
+ \omega_n(\alpha),
\end{equation}
where $ \sqrt{n}( H_n(\xi^+_\alpha)-H_F(\xi^+_\alpha)
) \stackrel{d}{\longrightarrow} \tau^+_\alpha$ and
$ \omega_n(\alpha) = \sqrt{n}N^{-1}(\alpha N - [\alpha N]) =
\mathrm{O}(n^{-1/2}) $
as \mbox{$ n \rightarrow\infty$}, but the positive term
$\sqrt{n}\delta_\alpha(\xi^+_\alpha,0) $ is unbounded.
Therefore, in this case, we shall apply the estimate (\ref{for14}) with $
\varepsilon
n^{-1} $ instead of $ \varepsilon, $
that is, $ P\{ (h_{nN_\alpha}-\xi^+_\alpha> \varepsilon n^{-1} )$ $
\cap\ (h_{nN_\alpha} >\xi^-_\alpha)\} $
$ \leq c_1 \exp\{ -c_2 n \delta^2_\alpha(\xi^+_\alpha,\varepsilon
n^{-1})\}.$
Since the distribution function $ H_F $ is continuous from the right at
the point $ \xi^+_\alpha$
it follows that $ \delta_\alpha(\xi^+_\alpha,0)  \leq\delta
_\alpha(\xi^+_\alpha,\varepsilon n^{-1}) $ for any small $
\varepsilon> 0 $
and sufficiently large~$ n.$ Hence, in the second case, $ \alpha<
H_F(\xi^+_\alpha) $ and (\ref{for14}) is replaced by the inequality
%
\begin{equation}\label{for16}
P\{ (h_{nN_\alpha}-\xi^+_\alpha> \varepsilon n^{-1} ) \cap
(h_{nN_\alpha} >\xi^-_\alpha)\}
\leq c_1 \exp\{ -c_2 n \delta^2_\alpha(\xi^+_\alpha,0)\},
\end{equation}
which provides the desired convergence $ \sqrt{n}N^{-1} \bar J_\alpha
\rightarrow0 $
in probability as $ n \rightarrow\infty.$ Thus, we have proven that $
\sqrt{n}N^{-1} \mathbb{L}_\alpha\rightarrow0 $
in probability as $ n \rightarrow\infty.$

\textbf{Estimating $ \mathbb{L}_\beta$}. Noting that $ I\{\xi^-_\beta
\geq h_{ni} \} = I\{i\leq N^-_\beta\} $,
we write
\begin{eqnarray*}
0 &\leq&- J_\beta =  I\{N_\beta< N^-_\beta\}\sum_{i=N_\beta
+1}^{N^-_\beta}(\xi^-_\beta- h_{ni})
I\{\xi^-_\beta\geq h_{ni}\} \nonumber\\
& \leq& I\{N_\beta< N^-_\beta\}( N^-_\beta- N_\beta)(\xi^-_\beta
- h_{nN_\beta})
I\{\xi^-_\beta\geq h_{nN_\beta}\}.\nonumber
\end{eqnarray*}
Here, by analogy with (\ref{for12}), we have
%
\begin{eqnarray}\label{for17}
&&P\{ (\xi^-_\beta- h_{nN_\beta} > \varepsilon) \cap( \xi
^-_\beta\geq h_{nN_\beta} )\}\leq  P\{ \xi^-_\beta - h_{nN_\beta} \geq\varepsilon\}
\nonumber\\
&&\quad =  P\{ H_n(\xi^-_\beta-\varepsilon) \geq N^{-1}N_\beta\}
\nonumber\\ 
&&\quad  =  P\{ H_n(\xi^-_\beta-\varepsilon) - H(\xi^-_\beta-\varepsilon
) \geq N^{-1}N_\beta- H(\xi^-_\beta-\varepsilon)\} \nonumber\\
&&\quad  \leq c_1 \exp\{ -c_2 n \theta^2_\beta(\xi^-_\beta,\varepsilon
)\}
\end{eqnarray}
with $ \theta_\beta(\xi^-_\beta,\varepsilon) = \beta- H_F( \xi
^-_\beta
-\varepsilon)$ and by analogy with (\ref{for15}),
%
\begin{equation}\label{for18}
\sqrt{n}N^{-1}(N^-_\beta- N_\beta) = \sqrt{n}\bigl(H_n(\xi^-_\beta
)-H_F(\xi^-_\beta)\bigr)+\sqrt{n}\theta_\beta(\xi^-_\beta,0)
+ \omega_n(\beta).
\end{equation}
Here, we need to consider two cases: $ \beta- H_F( \xi^-_\beta- ) =
0 $
and $ \beta- H_F( \xi^-_\beta-)> 0.$ In the first case, we apply
the inequality (\ref{for17}) with sufficiently small $ \varepsilon> 0 .$ In the
second case, we use (\ref{for17}) again, but with parameter $ \varepsilon n^{-1}$, as in
(\ref{for16}), to get
%
\begin{equation}\label{for19}
P\{ (\xi^-_\beta- h_{nN_\beta} > \varepsilon n^{-1} ) \cap( \xi
^-_\beta\geq h_{nN_\beta} )\}
\leq c_1 \exp\{ -c_2 n \theta^2_\beta(\xi^-_\beta-,0)\}
\end{equation}
since the distribution function $ H_F $ has a limit from the left at the
point $ \xi^-_\beta$ and $H_F(\xi_{\beta}^{-}-) \geq H_F(\xi
_{\beta}^{-} -
\varepsilon n^{-1})$.
In this result, we have $ \sqrt{n}N^{-1} J_\beta\rightarrow0 $
in probability as $ n \rightarrow\infty.$

Finally, we consider $ \bar J_\beta.$ Since
$ I\{h_{ni} > \xi^\pm_\beta\} = I\{i > N^\pm_\beta\} $, we write
\begin{eqnarray*}
\bar J_\beta& = & I\{N^-_\beta< N_\beta\}\sum_{i=N^-_\beta
+1}^{N_\beta}(h_{ni}-\xi^+_\beta)I\{h_{ni} > \xi^-_\beta\}
\nonumber\\[-1pt]
& = & - I\{N^-_\beta< N_\beta\}\sum_{i=N^-_\beta+1}^{N_\beta}(
\xi^+_\beta- h_{ni})I\{\xi^-_\beta< h_{ni} < \xi^+_\beta\} I\{
N^-_\beta< i \leq\dot N^+_\beta\}\nonumber\\[-1pt]
& & {}+ I\{N^-_\beta< N_\beta\}\sum_{i=N^+_\beta+1}^{N_\beta}(h_{ni} -
\xi^+_\beta)I\{h_{ni} > \xi^+_\beta\} I\{ i > \dot N^+_\beta\}
\nonumber\\[-1pt]
& = & -\bar J^-_\beta+ \bar J^+_\beta.
\end{eqnarray*}
If $ \xi^-_\beta= \xi^+_\beta$, then $ \bar J^-_\beta= 0 $ a.s.
Now, assume that $ \xi^-_\beta\neq\xi^+_\beta.$
In this case, $ H_F(\xi^-_\beta)= \beta= H_F(\xi^+_\beta-)
$ and
we have
\begin{eqnarray*}
0 &\leq&\bar J^-_\beta \leq I\{N^-_\beta< \dot N^+_\beta\}(\xi
^+_\beta- \xi^-_\beta)\sum_{i=N^-_\beta+1}^{\dot N^+_\beta}
I\{\xi^-_\beta< h_{ni} < \xi^+_\beta\} \\
& \leq& \Delta\xi_\beta \sum_{i=1}^{N}
I\{\xi^-_\beta< h_i < \xi^+_\beta\} = 0 \qquad    \mbox{a.s.}
\end{eqnarray*}
since $ EI\{\xi^-_\beta< h_i < \xi^+_\beta\} = H_F(\xi^+_\beta-) -
H_F(\xi^-_\beta) = 0.$ Hence, we always have
$ \bar J^-_\beta= 0 $ a.s. To estimate $ \bar J^+_\beta$,
we write
\[
0 \leq\bar J^+_\beta \leq I\{N^-_\beta< N_\beta\}(N^+_\beta-
N_\beta) (h_{ni} -
\xi^+_\beta)I\{h_{ni} > \xi^+_\beta\}
\]
and apply the estimates (\ref{for13})--(\ref{for16}) with $ \beta$ instead of $ \alpha.$
We have
$ \sqrt{n}N^{-1} \bar J_\beta\rightarrow0 $
in probability as $ n \rightarrow\infty$ and hence $ \sqrt{n}N^{-1}
\mathbb{L}_\beta\rightarrow0 $
in probability as $ n \rightarrow\infty.$

This proves Lemma \ref{lem2.2}.
\end{pf}

\begin{pf*}{Proof of Theorem \ref{teo1.1}}
Let $ U(g) $ be a $ U $-statistic of the form (\ref{for1}) with the kernel
\begin{eqnarray*}
g(x_1,\ldots,x_m) & = & \bigl[ I\{h(x_1,\ldots,x_m) \leq\xi
^-_\beta\}\bigl(h(x_1,\ldots,x_m) - \xi^-_\beta\bigr)+\beta\xi^-_\beta\bigr]
\nonumber\\
&&{} -  \bigl[ I\{ h(x_1,\ldots,x_m) < \xi^+_\alpha\}\bigl(h(x_1,\ldots,x_m) -
\xi^+_\alpha\bigr)+\alpha\xi^+_\alpha\bigr] . \nonumber
\end{eqnarray*}
We see that
\begin{eqnarray*}
U(g) & = & N^{-1} \sum_{i = 1}^{N} I\{\xi^+_\alpha\leq h_i
\leq\xi^-_\beta\}h_i  +  \xi^+_\alpha\bigl(H_n(\xi^+_\alpha-) - \alpha\bigr)
- \xi^-_\beta\bigl(H_n(\xi^-_\beta) - \beta\bigr) .\nonumber
\end{eqnarray*}
It is not difficult to verify for this function that
$ Eg(X_1,\ldots,X_m) = \theta$ and
$ g(x) = Eg(x,X_2,\ldots, X_m) - \theta,   x \in{\mathsf{X}}$; in
addition, $ Eg^2(X) > 0 $, by the condition of the theorem.
Hence, the kernel $ g $ is non-degenerate and, by the central limit
theorem for $ U $-statistics with such bounded kernels, we have the
weak convergence
$\tau_{ng}:= m^{-1}\sqrt{n}  ( U(\mathrm{g}) - \theta) \stackrel
{d}{\longrightarrow} \tau_g $ as $ n \rightarrow\infty $
(see, for example, Borovskikh \cite{Bo}).
By the same central limit theorem, we have
\[
\tau_{n\alpha}:= m^{-1}\sqrt{n}  \bigl( H_n(\xi^+_\alpha-)- H(\xi
^+_\alpha-) \bigr) \stackrel{d}{\longrightarrow} \tau_{\alpha}
\]
and
\[
\tau_{n\beta}:= m^{-1}\sqrt{n}  \bigl( H_n(\xi^-_\beta)- H(\xi
^-_\beta) \bigr) \stackrel{d}{\longrightarrow} \tau_{\beta}
\]
as $ n \rightarrow\infty.$ Under the conditions of the theorem, we
have
$ E|I\{\dot N^+_\alpha-N_\alpha> 0\}-I\{\tau_\alpha> 0\}|
\rightarrow0 $ if $ \Delta\xi_\alpha\neq0 $
(in this case, $ H( \xi^+_\alpha-)=\alpha$) and
$ E|I\{N^-_\beta- N_\beta< 0\}-I\{\tau_\beta< 0\}| \rightarrow0 $
if $ \Delta\xi_\beta\neq0 $
(in this case, $H(\xi^-_\beta)=\beta$).
Further, it is easy to prove that the covariances
$ \operatorname{Cov}(\tau_{n\ast},\tau_{n\star})\rightarrow\operatorname{Cov}(\tau_{\ast},\tau_{\star}) $ as $ n \rightarrow\infty $,
where $ \ast, \star= \alpha, g, \beta.$ Now, apply Lemma \ref{lem2.2} to
complete the
proof of Theorem \ref{teo1.1}.
\end{pf*}

\begin{Lem}\label{lem2.3}
The following representation holds:
\[
L_{\alpha\beta} = N^{-1}\sum_{i = \bar{N}_{\alpha}}^{\bar
{N}_{\beta} -1}  h_{ni}.
\]
\end{Lem}

\begin{pf}
By definition,
we can write
\begin{eqnarray*}
L_{\alpha\beta} & = & \int_{R}I\{h_\alpha \leq x < h_\beta\} x
\,\mathrm{d}H_n(x) \\
&= & \frac{1}{N} \sum_{i = 1}^{N} I\{h_\alpha\leq h_{ni} < h_\beta
\} h_{ni}\\
&=& \frac{1}{N} \sum_{i=1}^N I\{ h_{ni} < h_{\beta} \} h_{ni} -
\frac{1}{N}
\sum_{i=1}^N I\{ h_{ni} < h_{\alpha} \} h_{ni} \\
&=& \frac{1}{N} \sum_{i=1}^{\bar{N}_{\beta} -1}  h_{ni} - \frac{1}{N}
\sum_{i=1}^{\bar{N}_{\alpha} -1}  h_{ni} \\
&=& \frac{1}{N} \sum_{i=\bar{N}_{\alpha}}^{\bar{N}_{\beta} -1}  h_{ni} .
\end{eqnarray*}
This proves Lemma \ref{lem2.3}.
\end{pf}

The proof of Corollary \ref{cor1.2} follows from Theorem \ref{teo1.1} and Lemma \ref{lem2.3}.

\printhistory

\end{document}